\documentclass[11pt,twoside, final]{amsart}
\copyrightinfo{0}{Iranian Mathematical Society}
\pagespan{1}{\pageref*{LastPage}}
\usepackage{etoolbox,lastpage}
\commby{}
\date{\scriptsize   Received: 13 September 2016 , Accepted: 23 March 2017.}
\usepackage{amsmath,amsthm,amscd,amsfonts,amssymb,enumerate}
\usepackage{graphicx}		
\usepackage{color}
\usepackage[colorlinks]{hyperref}
 \usepackage{euscript}
 
 \newcommand{\cA}{\EuScript{A}}

\newcommand{\cC}{\EuScript{C}}

\newcommand{\cF}{\EuScript{F}}

\newcommand{\cL}{\EuScript{L}}
\newcommand{\cM}{\EuScript{M}}
\newcommand{\cN}{\EuScript{N}}

\newcommand{\cR}{\EuScript{R}}
\newcommand{\cS}{\EuScript{S}}

\newcommand{\cV}{\EuScript{V}}

\newtheorem{theorem}{Theorem}[section]

\newtheorem{lemma}[theorem]{Lemma}

\theoremstyle{definition}

\theoremstyle{remark}

\numberwithin{equation}{section}
 
\begin{document}

 
\title[Theorems of Burnside and  Wedderburn]{Theorems of Burnside and Wedderburn revisited}
 
\author[B.R. Yahaghi]{Bamdad R. Yahaghi}
\address{
Department of Mathematics, Faculty of Sciences, Golestan University, Gorgan 19395-5746, Iran}
\email{ bamdad5@hotmail.com, bamdad@bamdadyahaghi.com}

%
 
\dedicatory{ Dedicated to the memory of my beloved parents, the light of my heavens and earth,\\
``for the flame that never dies will always remain in the heart..."} 
 
 \maketitle
%

\begin{abstract}
We approach celebrated theorems of Burnside and Wedderburn via simultaneous triangularization. First, for a general field $F$, we prove that $M_n(F)$ is the only irrreducible subalgebra of triangularizable matrices in $M_n(F)$ provided such a subalgebra exists. This provides a slight generalization of a well-known theorem of Burnside.  Next, for a given $n > 1$, we characterize all fields $F$ such that Burnside's Theorem holds in $M_n(F)$, i.e., $M_n(F)$ is the only irreducible subalgebra of itself. In fact, for a subfield $F$ of the center of a division ring $D$, our simple proof of the aforementioned extension of Burnside's Theorem can be adjusted to establish a Burnside type theorem for irreducible $F$-algebras of triangularizable matrices in $M_n(D)$ with inner eigenvalues in $F$, namely such subalgebras of $M_n(D)$ are similar to $M_n(F)$. We use Burnside's theorem to present a simple proof of a theorem due to Wedderburn. Then, we use our Burnside type theorem to prove an extension of Wedderburn's Theorem as follows: A subalgebra of a semi-simple left Artinian $F$-algebra  is nilpotent iff the algebra, as a vector space over the field $F$, is spanned by its nilpotent members and that the minimal polynomials of all of its members split into linear factors over $F$. We conclude with an application of Wedderburn's Theorem.
\end{abstract}
 
\section{\bf Introduction}

In 1905, W. Burnside proved  a  theorem, which is now a standard result, asserting that a group of $n \times n$ complex matrices is irreducible iff it contains a vector space basis for $M_n(\mathbb{C})$, equivalently, its linear span is $M_n(\mathbb{C})$, see \cite[Theorem on p.433]{B}. Later on, this result of Burnside was generalized by Frobenius and Schur in \cite{FS}, eventually giving rise to a modern form of Burnside's result from an algebra-theoretic perspective, which is: $M_n(\mathbb{C})$ is the only irreducible subalgebra of itself. Nowadays the following is known as Burnside's Theorem: for an algebraically closed field $F$, $M_n(F)$ is the only irreducible subalgebra of itself. Clearly, one can naturally state  a counterpart of Burnside's Theorem for irreducible algebras of linear transformations on finite-dimensional vector spaces over algebraically closed fields  \cite[Theorem 1.2.2]{RR}. Many standard triangularization results for collections of matrices with entries from algebraically closed fields (resp. linear transformations on finite-dimensional vector spaces over algebraically closed fields) can be proved using Burnside's Theorem, e.g., Theorems of McCoy, Kaplansky, Kolchin and the Block Triangularization Theorem to name a few; see  \cite[Chapters 1-4]{RR}. Also, Burnside's Theorem has applications  in linear groups and hence in the representation theory of groups. For an exposition of such applications of Burnside's Theorem, see \cite[Part II.2]{K} and \cite[\S9]{La}. A well-known theorem of J.H.M. Wedderburn states that a finite-dimensional algebra over a general field is nilpotent if the algebra, as a vector space over the field, is spanned by its nilpotent members;  see \cite[Theorem 1]{W}, \cite[Proposition 4.6]{P}, and \cite[Proposition 2.6.32]{R}.  Standard proofs of this theorem of Wedderburn, of which  the author is aware,  make  use of noncommutative ring theory.

In this note, we approach the aforementioned celebrated theorems of Burnside and Wedderburn via simultaneous triangularization and provide simple proofs for the theorems. Motivated by the algebra-theoretic version of Burnside's Theorem, for a general field $F$, it is shown that $M_n(F)$ is the only irreducible algebra  of triangularizable matrices in $M_n(F)$  provided such a subalgebra exists  (see below for definitions). We refer the reader to \cite{B}, \cite{FS}, \cite{HR}, \cite{K}, \cite{L}, \cite{LZ}, \cite{LR}, \cite{RRY}, \cite{RRo}, \cite{RR}, \cite{RY}, \cite{Ro},  \cite{Y1}, \cite{Y2}, and \cite{Y3} for proofs and extensions of Burnside's Theorem. Then, for a given $n > 1$, we characterize all fields $F$ such that Burnside's Theorem holds in $M_n(F)$. Indeed, our proof of the aforementioned extension of Burnside's Theorem  has at least two advantages. (a) It essentially proves that an algebra of matrices in $M_n(F)$ with $F$ a field over which there exists an irreducible quadratic polynomial, e.g., finite extensions of prime fields or the real field, is triangularizable iff it consists of triangularizable matrices. (b) Our proof has the advantage that it can be adjusted to establish a Burnside type theorem for $F$-algebras of triangularizable matrices with entries from a general division ring having inner-eigenvalues in the subfield $F$ of the center of the division ring; see Theorem \ref{2.3} below. With this result at disposal, one can develop simultaneous triangularization theory for collections of matrices with entries from general division rings having inner-eigenvalues in the centers of the division rings. Next, with  Burnside's theorem and the standard techniques of simultaneous triangularization at our disposal, we present a simple proof of  Wedderburn's Theorem. We then use our Burnside type theorem, namely Theorem \ref{2.3}, to prove an extension of Wedderburn's Theorem as follows: An $F$-algebra (not necessarily finite-dimensional) that is embedded in a semi-simple left Artinian $F$-algebra is nilpotent iff the algebra, as a vector space over the field $F$, is spanned by its nilpotent members  and that the minimal polynomials of all of its members split into linear factors over $F$; see Theorem \ref{2.5} and its consequences as pointed out in the remark following the theorem. While this generalization of Wedderburn's Theorem  is essentially a consequence of the proof of our extension of Burnside's theorem, it is certainly beyond the scope of the standard version of Burnside's Theorem. This can be thought of as a third advantage of the proof of our extension of Burnside's theorem. We conclude with an application of Wedderburn's Theorem.

Let us start setting the stage by establishing some definitions
and standard notation. Throughout the paper, $D$ denotes a
general division ring and $F$ stands for a
general field or a subfield of the center of $D$. As is usual,
we use the symbol ${\rm Z}(D)$ to denote the center of $D$. We
view the members of $M_n(D)$ as linear transformations acting on
the left of $D^n$, where $D^n$ is the right vector space of $n
\times 1$ column vectors.  The subspaces $\{0\}$ and $D^n$ are
called the trivial subspaces of $D^n$.

A subspace \(\cM\) is  {\it invariant} for a collection \(\cF\) of
matrices if $T{\cM \subseteq \cM} $  for all $ T \in \cF$. A
collection \(\cF\) of matrices in $M_n(D)$ is called {\it
reducible} if $\cF = \{0\}$ or it has a nontrivial invariant
subspace. The collection \(\cF\) is called irreducible if it is not reducible.

A collection $\cF$ of matrices in $M_n(D)$ is called {\it
simultaneously triangularizable} or simply {\it triangularizable}
if there exists  a basis for the vector space $D^n$, called a triangularizing basis for $\mathcal{F}$, relative to
which all matrices in the family are upper triangular.
Equivalently, there exists an invertible matrix $S$ over $D$ such
that each member of $S^{-1} \cF S$ is upper triangular, or there exists a maximal chain of invariant subspaces of $D^n$, called a triangularizing chain of subspaces for $\mathcal{F}$, each of which is invariant under $\mathcal{F}$.

For a triangularizable matrix $A \in M_n(D)$, we say that $\lambda \in D$ is an {\it inner
eigenvalue of $A$ relative to a member $\cM$ of a triangularizing
chain} $\cC$ for $A$ if there exists a column vector $x \in  \cM \setminus \cM_{-}$
such that $A x - x \lambda\in \cM_{-}$, where  $\cM_{-}$ is  the
predecessor of $\cM$ in $\cC$ (note that $\dim \cM/\cM_{-} =1$).
If $D=F$, then inner eigenvalues of a triangularizable matrix $A \in M_n(F)$ relative to the
members of a triangularizing chain $\cC$ for $A$ are the
eigenvalues of $A$. Inner eigenvalues of a triangularizable matrix $A \in M_n(D)$ relative
to a member $\cM$ of a triangularizing chain for $A$ are also
known as {\it diagonal coefficients relative to} $\cM$. If a triangularizable matrix $A \in
M_n(D)$ has inner eigenvalues in $F$ relative to the members of a
triangularizing chain for $A$, then it is easy to verify that the
inner eigenvalues of $A$ relative to the members of any other
triangularizing chain for $A$ are in $F$. It is also easily verified that $ A \in M_n(D)$ is triangularizable and has inner eigenvalues in $F$ iff $A$ has a minimal polynomial $ f \in F[X]$ which splits into linear factors over $F$.

By an {\it $F$-algebra} $\cA$ in  $ M_n(D)$, we mean a subring of
$ M_n(D)$ that is closed under scalar multiplication by the
elements of the subfield $F$. For a semigroup  $\cS$ in $ M_n(D)$,
we use ${\rm Alg}_F(\cS)$ to denote the $F$-algebra generated by
$\cS$. By ${\rm Alg}(\cS)$ we simply mean ${\rm Alg}_Z(\cS)$, where
$Z$ denotes the center of $D$. A matrix $A \in M_n(D)$ is called
{\it $F$-algebraic} if it is algebraic over the subfield $F$.

For a given field $F$ and $k \in \mathbb{N}$ with $k > 1$, we say that
$F$ is {\it $k$-closed} if every polynomial of degree $k$ over $F$
is reducible over $F$. It is plain that a field $F$ is
algebraically closed if and only if $F$ is $k$-closed for all $k
\in \mathbb{N}$ with $k > 1$. It can be shown that finite extensions of
prime fields, e.g., finite fields, are not $k$-closed for all $k
\in \mathbb{N}$ with $k > 1$.

Let $D, E$ be division rings. The division ring $E$ is called an
extension of $D$ if $D \subset E$ and ${\rm Z}(D) \subset  {\rm
Z}(E)$. A collection \(\cF\) of matrices in \(M_n(D)\) is called
{\it absolutely irreducible} if it is irreducible over all
extensions of $D$. If $D=F$, then, by Burnside's Theorem, the
collection \(\cF\) is absolutely irreducible if and only if \({\rm
Alg}(\cF)= M_n(F)\).  It is plain that an absolutely irreducible family of matrices in \( M_n(F)\) is irreducible and its commutant consists of scalars. The converse of this easily seen fact follows from the Wedderburn-Artin type theorem below, namely Theorem \ref{1.1}. Recall that the {\it commutant} of a collection $\mathcal{F}$ of matrices, denoted by $\mathcal{F}'$, is the set of all matrices that commute with every element of the collection   $\mathcal{F}$. A collection of matrices in $M_n(D)$ is said to be {\it defined over a subfield} $F$ of the division ring $D$ if the collection is similar to a collection of matrices in $M_n(F)$.

\bigskip

In what follows, we will make use of the following  Wedderburn-Artin type theorem which was proved in \cite{Y2} (see \cite[Theorem 2.2]{Y2}).

\bigskip

\begin{theorem} \label{1.1}
Let $n \in \mathbb{N}$, $D$ be a division ring, $F$ a subfield of its center,
and $\cA$ an irreducible $F$-algebra of $F$-algebraic  matrices in
$ M_n(D)$. Let $r \in \mathbb{N}$  be the smallest nonzero rank present
in $\cA$. Then, the integer $r$ divides $n$ and after a similarity
$\cA = M_{n/r}(\Delta)$, where  $\Delta$ is an irreducible division
$F$-algebra of $F$-algebraic matrices in $ M_r(D)$. In
particular, after a similarity, $\cA =M_n(\Delta_1)$, where $\Delta_1$ is
an $F$-algebraic subdivision ring of $D$, if and only if $r=1$. 
\end{theorem}

\bigskip

For s subset $\cC$ of a vector space $\cV$ over a division ring $D$, as is usual, we use $\langle \cC \rangle$ to denote the space spanned by the subset $\cC$.  We need the following useful lemma in our proof of  Theorem \ref{2.1} below (see  \cite[Lemma 2.1.12]{RR} and \cite[Lemma 4.2.1]{Y3}).

\bigskip

\begin{lemma}\label{1.2}
{\it Let $\cV$ be a finite-dimensional vector space over a division ring $D$, $\cS$ a semigroup in $ \cL(\cV)$, and  $T$ a nonzero linear transformation in  $ \cL(\cV)$. If  $\cS$ is irreducible, then so is $T\cS|_{\cR}$, where  $\cR= T\cV$ is the range of $T$.}
\end{lemma}
\bigskip

\begin{proof} If $\dim \cV = 1$, then the assertion trivially holds. So we may assume, with no loss of generality, that $\dim \cV >1$. There are now two cases to consider.

\bigskip

(a) ${\rm rank}(T)=1$.

\bigskip

To prove the assertion by contradiction suppose $T\cS|_{\cR}$ is reducible. Since $\dim \cR=1$ in this case, it follows from  definition that $T\cS|_{\cR} = \{0\}$. Therefore, $T\cS T= \{0\}$. Pick a nonzero $x \in \cV$ such that $Tx \not = 0$. Now either $\cS Tx = \{0\}$ in which case $\langle Tx \rangle$ is a nontrivial invariant subspace for $\cS$, or else
$$\langle \cS Tx \rangle  = \{ \sum_{i=1}^{k} S_i Tx_i: k \in \mathbb N, S_i \in \cS, x_i \in \langle x \rangle \ (1 \leq i \leq k)\}$$ is a nontrivial invariant subspace for $\cS$, because $T\cS T= \{0\}$ and  $\cS$ is a semigroup. This contradicts the hypothesis that $\cS$ is irreducible.

\bigskip

(b) ${\rm rank}(T)>1$.

\bigskip

To prove that $T\cS|_{\cR}$ is irreducible we use contradiction. Suppose that $T\cS|_{\cR}$ is reducible. So there exists a nontrivial subspace $\cM$ of $\cR= T\cV$ such that  $T\cS\cM \subseteq \cM $. Choose a nonzero $ x \in \cM$ and note that $T\cS\langle x \rangle \subseteq  \cM $. The subspace  $$\langle \cS x \rangle = \{ \sum_{i=1}^{k} S_i x_i: k \in \mathbb N, S_i \in \cS, x_i \in \langle x \rangle \ (1 \leq i \leq k)\}$$ is an invariant subspace of $\cS$. Furthermore, it is proper, for $T\cS \langle x \rangle \subseteq  \cM \subset \cR$.  If $\cS x= 0 $, then $\langle x\rangle$ is a nontrivial invariant subspace for $\cS$, otherwise $\langle \cS x\rangle$ would be a nontrivial invariant subspace for $\cS$. So in any event we conclude that $\cS$ is reducible, a contradiction.
\end{proof}

\bigskip

 A subspace is called {\it hyperinvariant} for a collections of matrices if it is invariant under the union of the collection and its commutant. A matrix $A$ in $M_n(F)$ ($n > 1$) is called {\it reducible} if $A$ as a linear transformation on $F^n$ is reducible, i.e., it has a nontrivial invariant subspace. The following exercise follows  from the Cayley-Hamilton Theorem: a matrix $A$ in $M_n(F)$ is irreducible if and only if the characteristic polynomial for $A$ is irreducible over $F$; if and only if every nonzero $x$ in $F^n$ is a cyclic vector for $A$, i.e., $\{x, Ax, \ldots, A^{n-1}x\}$ spans $F^n$. Furthermore,  $$\{A\}'=F[A]=\{f(A): f \in F[X] \ {\rm with} \ \deg(f)\leq n-1\},$$ where  $\{A\}'$ denotes the commutant of $A$ (see \cite[Lemma 1.1.1]{Y3}). In contrast to this exercise, we have the following, which is \cite[Lemma 2.2.20]{Y3}.

\bigskip

\begin{lemma}  \label{1.3}
{\it Let $F$ be a field and $n >1$. A matrix $A$ in $M_n(F)$ has no nontrivial hyperinvariant subspace if and only if the minimal polynomial for $A$ is irreducible over $F$. Furthermore, after a similarity, $\{A\}'=M_{\frac{n}{r}}(F[C])$, where $r= \deg({\rm m}_A)$ (divides $n$) and $C= {\rm C(m}_A)$ in  $M_r(F)$ denotes the companion matrix of the minimal polynomial of $A$.}
\end{lemma}

\bigskip

\begin{proof}  The ``only if'' part of the assertion is easy.  To prove it by contradiction suppose ${\rm m}_A$, the minimal polynomial of $A$, is reducible over $F$.  So there exists a polynomial $f \in F[X]$ different from ${\rm m}_A$ that divides ${\rm m}_A$. It is now easily seen that $\ker (f(A))$ is a nontrivial hyperinvariant subspace for $A$, a contradiction. To see the  ``if'' part, suppose that ${\rm m}_A$ is irreducible over $F$. From the Rational Canonical Form Theorem (see Theorem VII.4.2(i) and Theorem VII.4.6(i) of \cite{H}), we see that $r= \deg({\rm m}_A)$ divides $n$ and that $A$ is similar to a direct sum of copies of the companion matrix of ${\rm m}_A$. More precisely, $A$ is similar to $C \oplus \cdots \oplus C \in M_{\frac{n}{r}}(F[C])$, where $C$ denotes the the companion matrix of ${\rm m}_A$. With no loss of generality, we may assume that  $A = C \oplus \cdots \oplus C \in M_{\frac{n}{r}}(F[C])$. In view of the aforementioned elementary exercise, a straightforward calculation shows that $\{A\}'=M_{\frac{n}{r}}(F[C])$. Since ${\rm m}_A$ is irreducible over $F$, we conclude that $F[C]$ is an irreducible algebra in $M_r(F)$ which, in turn, implies the irreducibility of the algebra  $\{A\}'=M_{\frac{n}{r}}(F[C])$ in $M_n(F)$. That is, the matrix $A$ has no nontrivial hyperinvariant subspace.
\end{proof}

\bigskip

\section{\bf Main Results}

\bigskip

First, we present a simple proof of the following result which is slightly stronger than Burnside's Theorem.

\bigskip

\begin{theorem} \label{2.1}
 Let   $n \in \mathbb{N}$, $F$ be a field, and $\mathcal A$ an
irreducible algebra of triangularizable matrices in $ M_n(F)$. Then, $\mathcal A= M_n(F)$. Therefore, the
field $F$ is $k$-closed for each $k= 2, \ldots, n$ whenever $n > 1$, and Burnside's Theorem holds in $M_n(F)$.
\end{theorem}

\bigskip

\begin{proof}   Pick $ T \in \mathcal A$ with rank$\ T= r$, where $r$ is the minimal nonzero rank present in  $\mathcal A$.
First, we show that $ r=1$. To see this, note that $T\mathcal A|_{TF^n}$ is an irreducible
division algebra of  triangularizable linear transformations on $TF^n$. This implies that $T\mathcal A|_{TF^n}= F I_{TF^n}$, from which we obtain $ r= \dim TF^n = 1$, for $T\mathcal A|_{TF^n}$ is irreducible by Lemma \ref{1.2}. Next, we show that
$P^{-1} \mathcal A P = M_n(F)$ for some $ P \in M_n(F)$, implying that $\mathcal A  = M_n(F)$. Since $ r=1$, we see that, after a similarity, $E_{11} \in \mathcal A$, where $E_{11}$ is the standard $n \times n$ matrix with $1$ in the $(1,1)$-place and zero elsewhere. To see this, choose a rank-one matrix $T \in \mathcal{A}$ and an $A \in \mathcal{A}$ such that $TA|_{TF^n}=  I_{TF^n}$. Then $TA$ is clearly a rank-one idempotent, and hence is similar to $E_{11}$ by a standard exercise in linear algebra.  We prove the assertion by induction on $n$. If $n=1$, we have nothing to prove. Suppose that the assertion holds for matrices of size less than $n$. Let $\mathcal A$ be an irreducible algebra of triangularizable matrices in $M_n(F)$. As $E_{11} \in \mathcal A$, we can write
$$ E_{11} \mathcal A E_{11} = F \oplus 0_{n-1} \subseteq \mathcal A, \ (I - E_{11}) \mathcal A (I- E_{11}) = 0_1 \oplus \mathcal A_{1}
\subseteq \mathcal A, $$
where $ E_{11} = I_1 \oplus 0_{n-1} \in \mathcal A$ and $\mathcal A_{1} \subseteq M_{n-1}(F)$ is an irreducible algebra of triangularizable matrices in $M_{n-1}(F)$. It follows from the induction hypothesis that  $\mathcal A_{1} = M_{n-1}(F)$. Therefore,
$$ E_{11} \mathcal A E_{11} = F \oplus 0_{n-1} \subseteq \mathcal A, \ (I - E_{11}) \mathcal A (I- E_{11}) = 0_1 \oplus M_{n-1}(F)
\subseteq \mathcal A. $$
As $\mathcal A$ is irreducible, we see that  $(I - E_{11}) \mathcal A E_{11} \not= 0$ and $ E_{11} \mathcal A (I- E_{11}) \not= 0$. Thus, there are $ 2 \leq i_0 \leq n$ and $ 2 \leq j_0 \leq n$ such that $E_{i_0i_0} \mathcal A E_{11} \not= 0$ and $ E_{11} \mathcal A E_{j_0j_0} \not= 0$, implying that $E_{i_01} \in \mathcal A$ and $E_{1j_0} \in \mathcal A$. This, in turn, yields $E_{1i_0} = E_{1j_0}E_{j_0i_0} \in \mathcal A$ and $E_{j_01} = E_{j_0i_0} E_{i_01} \in \mathcal A$.  It thus follows that $E_{i1} = E_{ii_0} E_{i_01} \in  \mathcal A $ and $E_{1j} = E_{1j_0} E_{j_0j} \in  \mathcal A $  for all $ 1 \leq i, j  \leq n$.  Consequently,
$ E_{ij}= E_{i1} E_{1j} \in  \mathcal A$ for all $ 1 \leq i, j \leq n$. Therefore, $\mathcal A= M_n(F)$, as desired.
\end{proof} 

\bigskip

\noindent {\bf Remark.}
A consequence of Theorem \ref{2.1} is the following: {\it Let $F$ be a field that is not $2$-closed and $\mathcal{A}$ an algebra of matrices in $M_n(F)$. Then $\mathcal{A}$ is triangularizable iff $\mathcal{A}$ consists of triangularizable matrices} (see \cite[Theorem 2.8]{Y2}). This result does not follow from Burnside's Theorem because the ground field $F$ in Burnside's Theorem is assumed to be algebraically closed.

\bigskip

The following theorem, for a given $n > 1$, characterizes all fields $F$ such that Burnside's Theorem holds in $M_n(F)$.

\bigskip

\begin{theorem}  \label{2.2} Let $F$ be a field and $n >1$. The following are equivalent.

{\rm (i)}  The only irreducible algebra in $M_n(F)$ is  $M_n(F)$.

{\rm (ii)}  Every irreducible family of matrices in $M_n(F)$ is absolutely irreducible.

{\rm (iii)}  The commutant of every irreducible family of matrices in $M_n(F)$ consists of scalars.

{\rm (iv)}  Every nonscalar matrix in $M_n(F)$ has a nontrivial hyperinvariant subspace.

{\rm (v)}  The field $F$ is  $k$-closed for all $k$ dividing $n$ with $k > 1$.
\end{theorem}

\bigskip

\begin{proof}  ``(i) $\Longrightarrow$ (ii)'' Obvious.

 ``(ii) $\Longrightarrow$ (iii)'' Let $\cF$ be an irreducible family of matrices in $M_n(F)$ and $\cA$ denote the algebra generated by $\cF$. It is plain that $\cF'= \cA'$. On the other hand, by the hypothesis, we must have $\cA= M_n(F)$. Therefore, $\cF'= \cA'= M_n(F)'= \{cI_n: c \in F\}$, proving the assertion.

``(iii) $\Longrightarrow$ (iv)'' Use contradiction. Suppose that the nonscalar matrix $A$ in  $M_n(F)$ has no nontrivial hyperinvariant subspace. Therefore, $\{A\}'$ must be an irreducible family, in fact algebra, of matrices in $M_n(F)$. Since $A \in (\{A\}')'$, it follows from the hypothesis that $A$ is scalar, a contradiction.

``(iv) $\Longrightarrow$ (v)'' Use contradiction again. Suppose that there exists an irreducible polynomial $f$ over $F$ such that deg$f=r >1 $ divides $n$. Let $C \in M_r(F)$ denote the companion matrix for $f$. Set $A= C \oplus \cdots \oplus C \in M_{\frac{n}{r}}(F[C])$. It is plain that $A$ is nonscalar and that m$_A = f$. Now irreducibility of m$_A$ contradicts the hypothesis in view of Lemma \ref{1.3}.

``(v) $\Longrightarrow$ (i)'' Let $\cA$ be an irreducible algebra in $M_n(F)$. Let $D=F$ and let $r$ and $\Delta$ be as in Theorem \ref{1.1}.   Since $\Delta$ is a division algebra in $M_r(F)$, the minimal polynomial ${\rm m}_A$ of every $A \in \Delta$ is irreducible over $F$, implying that ${\rm deg(m}_A{\rm )}$ divides $r$ by the Rational Canonical Form Theorem (see Theorem VII.4.2(i) and Theorem VII.4.6(i) of \cite{H}).  Therefore, $r=1$ and $\Delta=F$ because $F$ is  $k$-closed for all $k$ dividing $n$ with $k > 1$. It now follows from Theorem \ref{1.1}, or from the proof of Theorem \ref{2.1}, that $\cA = M_n(F)$, as desired.
\end{proof}

\bigskip

\noindent {\bf Remark.}  Since the field of real numbers is $k$-closed whenever $k$ is an odd number greater than $1$, a quick consequence of the preceding theorem is what we can call {\it Burnside's Theorem for real matrices}:
 {\it  Let $n \in \mathbb N$ be an odd number. Then the only irreducible algebra of matrices in $M_n(\mathbb{R})$ is $M_n(\mathbb{R})$.}

\bigskip

In fact, mimicking the proof of Theorem \ref{2.1} above one can prove the following result, which can be thought of as an extension of Burnside's Theorem to irreducible $F$-algebras of triangularizable  matrices in $ M_n(D)$ with inner-eigenvalues in a subfield $F$ of the center of $D$.  As pointed out in \cite{Y2}, the following is a quick consequence of Theorem \ref{1.1}.  However, here we provide a detailed proof of this result for the sake of completeness on the one hand, and on the other hand to show the extent to which the proof of Theorem \ref{1.1} is simplified in this special case. For some consequences of this theorem, see \cite{Y2}.

\bigskip

\begin{theorem} \label{2.3}
Let $n \in \mathbb{N}$, $D$ be a division ring, $F$ a subfield of its center,
and $\mathcal A$ an irreducible $F$-algebra of triangularizable  matrices in
$ M_n(D)$ with inner-eigenvalues in $F$. Then, after a similarity,
$\mathcal A = M_{n}(F)$. In
particular, $\mathcal A$ is defined over $F$, $\mathcal A$ is
absolutely irreducible, and  the
field $F$ is $k$-closed for each $k= 2, \ldots, n$ whenever $n > 1$. 
\end{theorem}

\bigskip

\begin{proof} Clearly, it suffices to show that $\mathcal A$ is similar to $ M_{n}(F)$.  To prove this, first, by an argument identical to that presented in the proof of Theorem \ref{2.1} but in the setting of vector spaces over division rings, we see that the minimal nonzero rank present in $\mathcal A$ is $1$. And hence, just as we saw in the proof of Theorem \ref{2.1}, $\mathcal{A}$ contains a rank-one idempotent, which is similar to $E_{11}$ by a standard exercise from linear algebra. Consequently, there exists an invertible matrix $ P \in M_n(D)$ such that $E_{11} \in P^{-1} \mathcal{A} P$. Next, again, as in the proof of Theorem \ref{2.1}, to prove the assertion, we proceed by induction on $n$, the size of matrices. If $n=1$, then the assertion trivially holds.  Assuming that the assertion holds for matrices of size less than $n$, we prove the assertion for matrices of size $n$. To this end, let $\mathcal A$ be as in the statement of the theorem. Clearly, in view of the hypothesis, we can write
$$ E_{11} P^{-1}\mathcal{A} P E_{11} = F \oplus 0_{n-1}  \subseteq P^{-1}\mathcal{A} P ,$$
$$ (I - E_{11}) P^{-1}\mathcal{A} P (I- E_{11}) = 0_1 \oplus \mathcal A_{1} \subseteq P^{-1}\mathcal{A} P , $$
where $ E_{11} = I_1 \oplus 0_{n-1} \in P^{-1}\mathcal{A} P$ and $\mathcal A_{1} \subseteq M_{n-1}(D)$ is an irreducible $F$-algebra of triangularizable matrices in $M_{n-1}(D)$ with inner-eigenvalues in $F$. It thus follows from the inductive  hypothesis that  there exists an invertible matrix $ Q_1 \in M_{n-1}(D)$ such that $ Q_1^{-1}\mathcal A_{1} Q_1 = M_{n-1}(F)$. Therefore, letting $ R := (I_1 \oplus Q_1) P \in M_n(D)$, we have $ E_{11} \in R^{-1} \mathcal A R$, and hence  we can write 
$$ E_{11} R^{-1} \mathcal A R E_{11} = F \oplus 0_{n-1} \subseteq R^{-1} \mathcal A R ,$$
$$ (I - E_{11}) R^{-1} \mathcal A R (I- E_{11}) = 0_1 \oplus M_{n-1}(F)
\subseteq R^{-1} \mathcal A R. $$
As $R^{-1} \mathcal A R$ is irreducible, we see that  
$$(I - E_{11})R^{-1} \mathcal A R E_{11} \not= 0 \ {\rm and} \  E_{11} R^{-1} \mathcal A R (I- E_{11}) \not= 0.$$
Thus, there are $ 2 \leq i_0 \leq n$ and $ 2 \leq j_0 \leq n$ such that $E_{i_0i_0} R^{-1} \mathcal A R E_{11} \not= 0$ and $ E_{11} R^{-1} \mathcal A R E_{j_0j_0} \not= 0$, implying that $ aE_{i_01} \in R^{-1} \mathcal A R$ and $bE_{1j_0} \in R^{-1} \mathcal A R$ for some nonzero $a, b \in D$. If necessary, by passing to $S^{-1}R^{-1} \mathcal A R S$ and replacing $RS$ with $R$ from the beginning, where  $ S:= {\rm diag}(b^{-1}, 1, \ldots, 1) \in M_n(D)$, we may assume  that $b =1$.   Now $E_{1i_0} = E_{1j_0} E_{j_0i_0}  \in R^{-1} \mathcal A R$, and hence $ aE_{11} =E_{1i_0}  aE_{i_01} \in R^{-1} \mathcal A R$, implying that $ a \in F$ because $\mathcal A$ has inner-eigenvalues in $F$. Thus, we might as well assume that $a=1$.   So far, we have 
$ E_{i_01}, E_{1i_0},  E_{1j_0} \in R^{-1} \mathcal A R$. We complete the proof by showing that $R^{-1} \mathcal A R= M_n(F)$.
To this end,  we have $E_{11} \in R^{-1} \mathcal A R$ and  $E_{i1} = E_{ii_0} E_{i_01} \in  R^{-1} \mathcal A R $ for all $ 2 \leq i \leq n$, implying that $E_{i1} \in  R^{-1} \mathcal A R $ for all $ 1 \leq i \leq n$. Likewise,  $E_{1j} \in  R^{-1} \mathcal A R $ for all $ 1 \leq j \leq n$. 
Consequently,
$ E_{ij}= E_{i1} E_{1j} \in  R^{-1} \mathcal A R$ for all $ 1 \leq i, j \leq n$. Therefore, $R^{-1} \mathcal A R \supseteq  M_n(F)$. Conversely, let $ A= (a_{ij} ) \in R^{-1} \mathcal A R \subseteq M_n(D)$ be arbitrary. Then $ a_{ij} E_{11} = E_{1i} A E_{j1} \in R^{-1} \mathcal A R$ for all $ 1 \leq i , j \leq n$, from which we see that  $ a_{ij} \in F$ for all $ 1 \leq i , j \leq n$, because  $\mathcal A$ has inner-eigenvalues in $F$. This yields $ R^{-1} \mathcal A R \subseteq M_n(F)$, completing the proof.  
\end{proof}

\bigskip

We now use Burnside's Theorem to prove the following theorem due to Wedderburn (see \cite[Theorem 1]{W}, \cite[Proposition 4.6]{P}, and \cite[Proposition 2.6.32]{R}). We recall that if $F$ is a field, every finite-dimensional $F$-algebra can be embedded in a unital algebra and that every $n$-dimensional unital $F$ algebra
$ \mathcal{A}$ can be embedded in $ M_n(F)$ via $ A \mapsto [L_A]_{\mathcal{B}}$, where $L_A: \mathcal{A} \longrightarrow \mathcal{A}$ is the left multiplication by $A$, i.e., $ L_A(B) = AB$ ($B \in \mathcal{A}$) and  $\mathcal{B}$ is a fixed vector space basis for  $ \mathcal{A}$. Therefore, every $n$-dimensional $F$-algebra can be embedded in $M_k(F)$, where $ k \in \{ n , n+1\}$. In fact, it can be shown that every $n$-dimensional $F$-algebra
$ \mathcal{A}$ can be embedded in $ M_n(F)$; see for instance \cite[Corollary 5.5.b]{P}.

\bigskip

\begin{theorem} \label{2.4}
Let $F$ be a field and $\mathcal A$ a finite-dimensional $F$-algebra which is spanned by nilpotents as a vector space over $F$. Then, the algebra $\mathcal A$  is nilpotent, more precisely, $\mathcal A^n = 0$, where $ n= \dim \mathcal A + 1$ and $\mathcal A^n:= \{ A_1 \cdots A_n: A_i \in \mathcal A, 1 \leq i \leq n\}$. In particular, $\mathcal A$ is a nilalgebra, i.e., consists of nilpotent elements.
\end{theorem}

\bigskip

\begin{proof}  As $ k= \dim \EuScript{A}$, we can embed $\EuScript{A}$ in $M_{k+1}(F)$ by the comment we made preceding the theorem. So without loss of generality, we can assume that $ A \subseteq M_n(F)$, where $ n = k+1$. As nilpotency does not depend on the ground field, we may assume, without loss of generality, that $F$ is algebraically closed. Note that the assertion trivially holds if $n=1$. So suppose that $n>1$.

First we show that the algebra $\EuScript{A}$ is triangularizable. Suppose that $\mathcal M \subset \cN$ are two invariant subspaces for $\mathcal A$ with $\dim \frac{\cN}{\cM} >1$. It suffices to show that  $\widehat{\mathcal{A}} $, the set of all quotient transformations $\hat A$ on $\frac{\cN}{\cM}$ ($A \in  \EuScript{A}$), is reducible. It easily follows from the hypothesis that every $\hat A \in   \widehat{\mathcal{A}} $ can be written as a linear combination of nilpotent elements of $\widehat{\mathcal{A}} $. Thus,  ${\rm tr}(\widehat{\mathcal{A}})=\{0\}$. Since the ground field $F$ is algebraically closed, from Burnside's Theorem, we see that $\widehat{\mathcal{A}} $ is reducible, as desired. Therefore, the algebra $\widehat{\mathcal{A}} $ is triangularizable. Thus, there exists a basis $\mathcal{B}$, for the  vector space $F^n$, with respect to which every $A \in \mathcal{A}$ has strictly upper triangular matrix representation, for $\mathcal A$ as a vector space over $F$ is spanned by its nilpotent elements. This shows that  $\mathcal{A}$  consists of nilpotents and hence  $\mathcal A^n = 0$ because every $A \in \mathcal{A}$ has strictly upper triangular matrix representation with respect to $\mathcal{B}$. This  completes the proof.
\end{proof} 

\bigskip

With Theorem \ref{2.3} at our disposal, we are ready to give an extension of Wedderburn's Theorem for certain algebras of matrices (not necessarily finite-dimensional)   over division rings that are spanned by their nilpotent elements, and consequently for certain subalgebras of semi-simple left Artinian $F$-algebras as explained in the remark that follows. 

\bigskip

\begin{theorem} \label{2.5}
Let $n \in \mathbb{N}$, $F$ be a subfield of the center of $D$, and $\mathcal A$ an $F$-algebra of triangularizable matrices in $M_n(D)$ with inner eigenvalues in $F$ that is spanned by its nilpotent elements as a vector space over $F$. Then, the algebra $\mathcal A$  is nilpotent. More precisely, $\mathcal A^n = 0$, where $\mathcal A^n:= \{ A_1 \cdots A_n: A_i \in \mathcal A, 1 \leq i \leq n\}$. In particular, $\mathcal A$ is a nilalgebra, i.e., consists of nilpotent elements.
\end{theorem}

\bigskip

\begin{proof} 
Clearly, with no loss of generality, we may assume that $ n > 1$. Let
$$0 = \mathcal{M}_0 <  \mathcal{M}_1 < \cdots < \mathcal{M}_k= D^n$$
be a maximal chain of invariant subspaces for the $F$-algebra $\mathcal A$. We prove the assertion by showing that $ \dim \frac{\mathcal{M}_i}{\mathcal{M}_{i-1}}= 1$ for all $ 1 \leq i \leq k$. Fix $ 1 \leq i \leq n$ and use $ \mathcal{A}_i$ to denote the set of quotient linear transformations induced by $\mathcal A$ on $\frac{\mathcal{M}_i}{\mathcal{M}_{i-1}}$. Clearly, $ \mathcal{A}_i$ is an irreducible  $F$-algebra of $F$-algebraic triangularizable linear transformations in $ \mathcal{L}( \frac{\mathcal{M}_i}{\mathcal{M}_{i-1}})$ with inner eigenvalues in $F$ that is spanned by its nilpotent elements. Fixing a basis for $ \frac{\mathcal{M}_i}{\mathcal{M}_{i-1}}$, we may think of $ \mathcal{A}_i$ as an irreducible $F$-algebra of triangularizable matrices in  $M_{n_i}(D)$ with inner eigenvalues in $F$, where $ n_i = \dim \frac{\mathcal{M}_i}{\mathcal{M}_{i-1}}$.  By Theorem \ref{2.3}, $\mathcal{A}_i \sim M_{n_i}(F)$. In particular, the $F$-algebra $\mathcal{A}_i $ is finite-dimensional and is spanned by its nilpotent elements. It thus follows from Theorem \ref{2.4} that $\mathcal{A}_i $ is nilpotent. This clearly implies $ n_i = 1$ and $\mathcal{A}_i = 0$ for all $ 1 \leq i \leq k$. Therefore $ k = n$ and the proof is complete.
\end{proof}

\bigskip

\noindent {\bf Remark.} Clearly, the theorem, i.e., the nilpotency of the $F$-algebra,  holds under the weaker hypothesis that the $F$-algebra $\mathcal A$ is embedded in a finite direct sum of matrix rings over division rings, where $F$ is a subfield of the intersection of the centers of the division rings, and that the minimal polynomials of all of its members split into linear factors over $F$; or equivalently, in view of the Wedderburn-Artin Theorem (\cite[Theorem IX.5.4]{H}),  the $F$-algebra $\mathcal A$ is embedded in a semi-simple left Aritinian $F$-algebra (as opposed to being embedded in a simple left Artinian $F$-algebra, which is the case in the statement of the theorem) and that the minimal polynomials of all of its members split into linear factors over $F$.

\bigskip

Let $\mathcal{A}$ be an $F$-algebra and $ \mathcal{C} \subseteq \mathcal{A}$. By $N_\sigma( \mathcal{C})$ we mean the set of all $ A \in \mathcal{A}$ that can be written as an $F$-linear combination of nilpotent elements from  $ \mathcal{C} $. By {\it a semigroup ideal}  $\mathcal J$ of a (multiplicative) semigroup $\mathcal S$, we mean a subset of $\mathcal S$ with the property that $ SJ \in \mathcal J$ and $ JS \in \mathcal J$ for all $J \in \mathcal J$ and $S \in \mathcal S$.
 We conclude with the following application of Wedderburn's Theorem.

\bigskip

\begin{theorem} \label{2.6}
Let $ n > 1$, $F$ be a field, $\mathcal S$ an irreducible semigroup in $M_n(F)$ and $\mathcal J$ a nonzero semigroup ideal of $\mathcal{S}$. Then

{\rm (i)}
$\Big\{A \in {\rm Alg}_{F}( \mathcal{S}\cup \{I \}) : \mathcal{J} A\mathcal{J} \subseteq N_\sigma (\mathcal{J} A \mathcal{J})  \Big\} = \{0\}.$

{\rm (ii)}
$\Big\{A \in {\rm Alg}_{F}( \mathcal{S}\cup \{I \}) : A\mathcal{J} \subseteq  N_\sigma (A \mathcal{J}) \Big\} = \{0\}.$

{\rm (iii)}
$\Big\{A \in {\rm Alg}_{F}( \mathcal{S}\cup \{I \}) : \mathcal{J} A \subseteq  N_\sigma (\mathcal{J} A ) \Big\} = \{0\}.$

{\rm (iv)}
$\Big\{A \in {\rm Alg}_{F}( \mathcal{S}\cup \{I \}) : (J_1AJ_2)^n = 0 \ \forall J_1 , J_2 \in \mathcal{J}\Big\} = \{0\}.$

{\rm (v)}
$\Big\{A \in {\rm Alg}_{F}( \mathcal{S}\cup \{I \}) : (AJ)^n = 0 \ \forall J \in \mathcal{J}\Big\} = \{0\}.$

{\rm (vi)}
 $\Big\{A \in {\rm Alg}_{F}( \mathcal{S}\cup \{I \}) : (JA)^n = 0 \ \forall J \in \mathcal{J}\Big\} = \{0\}.$
\end{theorem}

\bigskip

\begin{proof}  Clearly (iv), (v), and (vi)  follow from (i), (ii), and (iii),  respectively. We prove (i); (ii) and (iii) can be proved in a similar fashion. Let
$A \in {\rm Alg}_{F}( \mathcal{S}\cup \{I \})$ be such that $ \mathcal{J} A\mathcal{J} \subseteq N_\sigma (\mathcal{J} A \mathcal{J})$. We prove that $ A = 0$.
It is plain that $  \mathcal{J} A\mathcal{J} \subseteq {\rm Alg}(\mathcal{J})$. The algebra ${\rm Alg}(\mathcal{J})$ is irreducible because the semigroup $\mathcal{S}$ is irreducible and $ \mathcal{J}$ is a nonzero semigroup ideal of $\mathcal{S}$. On the other hand, it follows from the hypothesis that ${\rm Alg} (\mathcal{J} A\mathcal{J})$ is spanned by its nilpotent elements. So by Theorem \ref{2.4},  ${\rm Alg} (\mathcal{J} A\mathcal{J})$  is nilpotent, and hence reducible. This implies ${\rm Alg} (\mathcal{J} A\mathcal{J})= 0$, for  ${\rm Alg} (\mathcal{J} A\mathcal{J})$ is an ideal of the irreducible algebra ${\rm Alg}(\mathcal{J})$. Consequently, $\mathcal{J} A\mathcal{J}= 0$, from which we easily see that $ A = 0$ because the nonzero semigroup ideal $\mathcal J$ is irreducible. This completes the proof. 
\end{proof}

\bigskip

\section*{\bf Acknowledgments}
This paper was written while the author was on sabbatical leave at the University of Waterloo. He would like to thank Department of Pure Mathematics of the University of Waterloo, and in particular Professors Laurent Marcoux and Heydar Radjavi, for their support. The author would like to thank the referee for his/her helpful comments.

\end{document}